\documentclass{amsart}

\usepackage[margin=1.5in]{geometry}
\usepackage{amsmath, graphicx, cite, enumerate, comment}
\usepackage{amssymb, amsthm, mathrsfs}

\newtheorem{theorem}{Theorem}[section]

\newtheorem{lemma}[theorem]{Lemma}

\theoremstyle{definition}

\theoremstyle{remark}
\newtheorem{remark}[theorem]{Remark}

\numberwithin{equation}{section}

\newcommand{\tr}{\textup{tr}}
\newcommand{\Rm}{\textup{Rm}}
\newcommand{\Ric}{\textup{Ric}}
\newcommand{\Div}{\textup{div}}

\newcommand{\Tr}{\textup{Tr}}

\newcommand{\CP}{\mathbb{CP}}

\newcommand{\spec}{\textup{spec}}

\begin{document}
\title{On 4-Dimensional $J$-Invariant Shrinking Ricci Solitons }

\author{Frederick Tsz-Ho Fong}
\address{Department of Mathematics, Brown University}
\email{fong@math.brown.edu}

\date{January 18, 2014}
\begin{abstract}
As of today, there are very few known complete shrinking Ricci solitons in dimension 4, and all examples discovered so far are K\"ahler and/or Einstein. In this note, we prove that all four dimensional $J$-invariant gradient shrinking Ricci solitons satisfy a differential form identity relating K\"ahlerity and Einstein-ness.
\end{abstract}
\maketitle
\section{Introduction}
A gradient shrinking Ricci soliton $g$ on a Riemannian manifold $M$ is a complete metric which satisfies
\begin{equation}\label{eq:SRS}
\Ric +\nabla^2 f = \frac{1}{2}g
\end{equation}
for some smooth function $f$, called the soliton potential. They are natural generalizations of positive Einstein metrics in a sense that the latter are shrinking Ricci solitons with constant potential. The study of shrinking Ricci soliton plays a central role in the Hamilton--Perelman theory of Ricci flow for many reasons: They are fixed-point solutions of the Ricci flow in the moduli space of Riemannian metrics modulo diffeomorphisms and rescalings. Also from variational point of view, they are critical metrics of the Perelman's $\mathcal{W}$-functional. More importantly, they arise as limits of dilations of Type I singularities of the Ricci flow. Therefore, in order to understand possible Ricci flow singularities, it is imparative to construct examples of shrinking Ricci solitons. Readers may refer to the survey article \cite{CaoSurv} for both historical and recent development of Ricci solitons in various aspects.

In dimension 4, known examples of complete shrinking Ricci solitons are very few as of today. All examples found so far are either Einstein and/or K\"ahler. In the non-Einstein case, there are only three examples known, namely the Cao--Koiso's shrinkers on $\CP^2 \# (-\CP^2)$ \cite{CaoSoliton, KoisoSoliton}, the Feldman--Ilmanen--Knopf's shrinker on the total space of $\mathcal{O}(-1)$-bundle over $\CP^1$ \cite{FIK}, and the Wang--Zhu's shrinker on $\CP^2 \# (-2\CP^2)$ \cite{WangZhu}. All of them are K\"ahler metrics. The first two examples are constructed using $U(2)$-rotationally symmetry, also known as the Calabi's ansatz. The ODE analysis for constructing these two examples reveal that the completeness condition imposes some very rigid restrictions on the parameters of the solutions of the ODE system so that only two $U(2)$-invariant examples can be found.

While examples found so far are K\"ahler and/or Einstein or direct products, it is not known so far whether all complete shrinking Ricci solitons must satisfy at least one of these conditions. While many classification results are obtained under various geometric assumptions (see \cite{CaoSurv} for a list of some), a complete classification of four dimensional shrinking Ricci solitons have not been carried out in full generality. The construction problem of a shrinking soliton which is non-K\"ahler, non-Einstein and non-product was also purposed in \cite{CaoSurv}. To the best of the author's knowledge, this construction problem is still open.

In search of new examples, it might be advantageous to establish some obstruction identities so that one might use them to rule out certain possibilities. In this note, we focus on shrinking Ricci solitons which are  complex surfaces (i.e. four real dimensional), and the metric $g$ and the Ricci tensor $\Ric$ are $J$-invariant.

A metric $g$ on a complex manifold $(M, J)$ with complex structure $J$ is said to be \textit{Hermitian} if it is $J$-invariant, i.e. $g(JX, JY) = g(X, Y)$ for any $X, Y \in TM$. For any Hermitian metric $g$, one can define the Hermitian form $\omega(X, Y) := g(X, JY)$. If in addition we have $d\omega = 0$, then the metric is said to be \textit{K\"ahler} and $\omega$ can be called the K\"ahler form. The main result of this note is the following:

\begin{theorem}\label{thm:main}
Let $(M^4, J)$ be a complex surface (not necessarily compact) and $g$ is a gradient shrinking Ricci soliton satisfying \eqref{eq:SRS}. Suppose $g$ is Hermitian and the Ricci tensor $\Ric$ is $J$-invariant, then we have the following differential form identity:
\begin{equation}\label{eq:main}
\Ric_0 (\nabla f, \cdot) \wedge d\omega = 0,
\end{equation}
where $\Ric_0 : = \Ric - \frac{R}{4}g$ is the trace-free part of the Ricci tensor $\Ric$.
\end{theorem}
\begin{remark}
The Einstein condition is equivalent to $\Ric_0 = 0$, and for Ricci solitons, it is also equivalent to $\nabla f = 0$. Therefore, the identity \eqref{eq:main} imposes an obstruction on these solitons in terms of Einstein-ness and K\"ahlerity.
\end{remark}
\begin{remark}
It is worthwhile to note that by the soliton equation \eqref{eq:SRS}, any two of the following conditions imply the third one:
\begin{itemize}
\item $g$ is Hermitian;
\item $\Ric$ is $J$-invariant;
\item $\nabla f$ is real holomorphic.
\end{itemize}
\end{remark}
\begin{remark}
The hypothese of Theorem \ref{thm:main} do not imply the metric $g$ must be K\"ahler, since there exist Einstein, Hermitian but non-K\"ahler, metrics on $\CP^2 \# (-\CP^2)$ and $\CP^2 \# (-2\CP^2)$, namely the Page metric \cite{Page} and the Chen-LeBrun-Weber's metric \cite{CLW} respectively.
\end{remark}
One key ingredient of the argument is to apply the Riemannian analogue of the Goldberg-Sachs Theorem \cite{AG} to show that the hypothesis of Theorem \ref{thm:main} implies that the self-dual Weyl tensor $W^+$ is degenerate (i.e. at every point at least two eigenvalues of $W^+$ are equal). Therefore, on the open dense subset of $M$ consisting of points at which $|\spec(W^+)|$ is locally constant, we are reduced to consider only the cases where $|\spec(W^+)| = 1$ or $2$.

In this note, we will first rewrite some geometric quantities in terms of the Hodge-star decomposition specific in four dimensions. The proof of Theorem \ref{thm:main} will be presented in Section 4.

\section{Riemann and Weyl curvature operators}
Let $(M^4, g)$ be an oriented Riemannian 4-manifold. The bundle of 2-forms over $M$ splits as
$$\wedge^2 M = \wedge^+ \oplus \wedge^-$$
where $\wedge^{\pm}$ is the $(\pm 1)$-eigenspace bundles of the Hodge star operator.

In local coordinates, the Weyl $(0,4)$-tensor can be expressed as
\begin{equation*}
W_{ijkl} = R_{ijkl} - \frac{1}{2} (g_{ik} R_{jl} + g_{jl} R_{ik} - g_{jk} R_{il} - g_{il} R_{jk}) + \frac{R}{6} (g_{ik} g_{jl} - g_{jk} g_{il}).
\end{equation*}
Here $R$ is the scalar curvature, $R_{ij}$'s are the components of the Ricci tensor $\Ric$ and $R_{ijkl}$'s are the components of the Riemann $(0,4)$-tensor. It is easy to verify that the Weyl tensor satisfies the relations
$$W_{ijkl} = W_{klij} = -W_{jikl}.$$
Therefore, the Weyl tensor can be considered as a self-adjoint endomorphism of $\wedge^2 M$ by decreeing
$$W(\omega)_{ij} = \frac{1}{2}W_{ijkl}\omega^{kl}$$
for any two-form $\omega$. The endomorphism $W$ commutes with the Hodge star operator $*$ on $\wedge^2 M$ and therefore $W$ perserves the decomposition $\wedge^{+} \oplus \wedge^{-}$. We denote $W^{\pm}$ to be the restriction of Weyl curvature operator $W$ on $\wedge^{\pm}$ and therefore we have $W = W^+ + W^-$. It is worthwhile to note that $W^\pm$'s are trace-free, i.e.
$$\Tr(W^\pm) = 0.$$

Given any local orthonormal frame $\{e_1, e_2, e_3, e_4\}$ on an open set of $M^4$, one can form a local frame $\{E_\alpha, F_\beta\}_{\alpha, \beta=1}^3$ respecting the Hodge-star decomposition in the following way:
\begin{align}\label{eq:EF}
\nonumber E_1 & = e_1 \wedge e_2 + e_3 \wedge e_4 & \qquad F_1 & = e_1 \wedge e_2 - e_3 \wedge e_4\\
E_2 & = e_1 \wedge e_3 + e_4 \wedge e_2 & \qquad F_2 & = e_1 \wedge e_3 - e_4 \wedge e_2\\
\nonumber E_3 & = e_3 \wedge e_2 + e_4 \wedge e_1 & \qquad F_3 & = e_3 \wedge e_2 - e_4 \wedge e_1
\end{align}
The $E_\alpha$'s span $\wedge^+$ and $F_\beta$'s span $\wedge^-$. The identification $\wedge^2_x M^4 \cong \mathfrak{so}(4)$ is given by, for instance,
\begin{equation*}
E_1 \sim \begin{bmatrix}
0&1&0&0\\
-1&0&0&0\\
0&0&0&1\\
0&0&-1&0
\end{bmatrix}.
\end{equation*}
The inner product between any pair of 2-forms $A, B$ is defined by the Killing form: $\langle A, B \rangle = \frac{1}{2} \Tr(AB)$. It can be easily verified that via this identification these local frame satisfy the following quaternionic relations:
$$E_1^2 = E_2^2 = E_3^2 = -I, \quad E_1 E_2 = E_3,$$
$$F_1 ^2 = F_2^2 = F_3^2 = -I, \quad F_1 F_2 = F_3.$$
Here the products and squares between the $E_\alpha$'s and $F_\beta$'s are simply the matrix products and squares. This observation will become handy for some computations later on.

Since the decomposition $\wedge^2 M^4 = \wedge^+ \oplus \wedge^-$ is invariant under parallel transport, one can express the covariant derivatives of $E_\alpha$'s in the following way:
\begin{equation}\label{eq:connection_E}
\nabla E_\alpha = \sum_{\beta=1}^3 \Omega_\alpha^\beta \otimes E_\beta.
\end{equation}
where $\{\Omega_\alpha^\beta\}_{\alpha, \beta=1}^3$'s are local connection 1-forms.

By the fact that $\langle E_\alpha, E_\beta \rangle = \delta_{\alpha\beta}$ and the compatibility of $\nabla$, we have
\begin{align}
\label{eq:switch_E}\Omega_\alpha^\alpha & = 0\\
\nonumber\Omega_\alpha^\beta & = -\Omega_\beta^\alpha.
\end{align}
Analogously, there exist 1-forms $\Gamma_\alpha^\beta$'s for the covariant derivatives of $F_\alpha$'s:
\begin{equation*}\label{eq:connection_F}
\nabla F_\alpha = \sum_{\beta=1}^3 \Gamma_\alpha^\beta \otimes F_\beta
\end{equation*}
and we have
\begin{align*}\label{eq:switch_F}
\Gamma_\alpha^\alpha & = 0\\
\nonumber\Gamma_\alpha^\beta & = - \Gamma_\beta^\alpha.
\end{align*}

From now on, we will use Greek letters $\alpha, \beta, \ldots \in \{1, 2, 3\}$ to denote the indices in $\wedge^2 M$ and Latin letters $i, j, k, \ldots \in \{1,2,3,4\}$ for the indices of $TM$. For convenience, the $\alpha$-indices are set to be modulo $3$, i.e. $E_{\alpha + 3} := E_\alpha$.

At a point $x \in M$, we diagonalize $W$ at $x$:
$$W = \frac{1}{2}\sum_{\alpha=1}^3 \lambda_\alpha E_\alpha \otimes E_\alpha + \frac{1}{2}\sum_{\beta=1}^3 \mu_\beta F_\beta \otimes F_\beta.$$
Here $\{E_\alpha\}_{\alpha=1}^3$ is an orthogonal basis of $\wedge^+$ at $x$ and $\{F_\beta\}_{\beta=1}^3$ is an orthogonal basis of $\wedge^-$. Each of $E_\alpha$'s and $F_\beta$'s has norm $\sqrt{2}$. By the traceless property of $W^\pm$, we have $$\sum_{\alpha=1}^3 \lambda_\alpha = \sum_{\beta=1}^3 \mu_\beta = 0.$$

Let $\mathcal{O}$ be the open and dense subset of $M$ consisting of points at which the number of distinct eigenvalues of $W^+$, denoted by $|\spec(W^+)|$, is locally constant. Then, $W^+$ is smoothly, locally diagonalizable on open subsets of $\mathcal{O}$ and the orthogonal basis $\{E_\alpha, F_\beta\}$ of $\wedge^2_x M^4$ extends smoothly to an orthogonal moving frame in an open neighborhood of $x$ inside $\mathcal{O}$. The eigenvalues $\lambda_\alpha$'s are also smoothly differentiable in this neighborhood. By \cite{AHS}, any orthogonal basis $\{E_\alpha, F_\beta\}$ of $\wedge^+ \oplus \wedge^-$ with norm $\sqrt{2}$ must be of the form \eqref{eq:EF} for some local orthonormal frame $\{e_i\}_{i=1}^4$ near $x$.

The Riemann curvature operator $\Rm$ acts on the local frame $\{E_\alpha, F_\beta\}$ in the following way. Recall that
$$R_{ijkl} = W_{ijkl} + \frac{1}{2}(g_{ik}R_{jl} + g_{jl}R_{ik} - g_{jk}R_{il} - g_{il}R_{jk}) - \frac{R}{6} (g_{ik}g_{jl} - g_{jk}g_{il}).$$
Fix an $\alpha$, one can verify:
\begin{equation*}
R_{ijkl} [E_\alpha]_{kl} = 2[W(E_\alpha)]_{ij} - \frac{R}{3} [E_\alpha]_{ij} + [E_\alpha \cdot \Ric + \Ric \cdot E_\alpha]_{ij}
\end{equation*}
where $R_{ijkl}$ denotes the component $\Rm(e_i, e_j, e_k, e_l)$ and $[E_\alpha]_{kl}$ is the $(k,l)$-th entry of the matrix in $\mathfrak{so}(4)$ identified with $E_\alpha$. Similarly, we have
$$R_{ijkl} [F_\alpha]_{kl} = 2[W(F_\alpha)]_{ij} - \frac{R}{3}[F_\alpha]_{ij} + [F_\alpha \cdot \Ric + \Ric \cdot F_\alpha]_{ij}.$$

The Cotton tensor $C_{ijk}$ is defined by
$$C_{ijk} = \nabla_k R_{ij} - \nabla_j R_{ik} - \frac{1}{6}(g_{ij}\nabla_k R - g_{ik}\nabla_j R).$$
From the second Bianchi identity the Cotton tensor is a constant multiple of the divergence of the Weyl curvature tensor $W$, precisely:
$$C_{ijk} = -2\Div(W)_{ijk} = -2\nabla_l W_{lijk}.$$
Let $(M^4, g, f)$ be a shrinking Ricci soliton satisfying \eqref{eq:SRS}. It is well-known (see e.g. \cite{ChowBook} for a proof) that any Ricci soliton satisfies:
\begin{equation}\label{eq:RicGradf}
2R_{ij}\nabla_j f = \nabla_i R.
\end{equation}
Combining \eqref{eq:SRS} and \eqref{eq:RicGradf}, the Cotton tensor for a Ricci soliton can be rewritten as
\begin{align*}
C_{ijk} & = -\nabla_k \nabla_i \nabla_j f + \nabla_j \nabla_i \nabla_k f - \frac{1}{3}(g_{ij}R_{kl}\nabla_l f - g_{ik} R_{jl} \nabla_l f)\\
& = -R_{pijk} \nabla_p f - \frac{1}{3}(g_{ij}R_{kl}\nabla_l f - g_{ik} R_{jl} \nabla_l f).
\end{align*}
Therefore, the Cotton tensor acts on the local frame $E_\alpha$ in the following way:
\begin{align*}
C_{ijk} [E_\alpha]_{jk} & = -R_{pijk} [E_\alpha]_{jk} \nabla_p f - \frac{1}{3} \left([E_\alpha]_{ik} R_{kl}\nabla_l f - [E_\alpha]_{ji} R_{jl} \nabla_l f\right)\\
& = 2\left([W(E_\alpha)]_{ip}-\frac{R}{6}[E_\alpha]_{ip} + \frac{1}{2} [E_\alpha \cdot \Ric + \Ric \cdot E_\alpha]_{ip}\right)\nabla_p f\\
& \quad - \frac{2}{3} [E_\alpha \cdot \Ric (\nabla f)]_i\\
& = 2[W(E_\alpha)]_{ip}\nabla_p f + \frac{1}{3}\left[E_\alpha \cdot \left(\Ric - R)(\nabla f)\right)\right]_i + [\Ric \cdot E_\alpha (\nabla f)]_i.
\end{align*}
Let $C^+$ be the restriction of $C: \wedge^1 \otimes \wedge^2$ on $\wedge^1 \otimes \wedge^+$. We then have
\begin{align}
\nonumber C^+ & = \frac{1}{2} \sum_{\alpha=1}^3 C^+(E_\alpha) \otimes E_\alpha\\
\nonumber & = \frac{1}{2} \sum_{\alpha=1}^3 \sum_{i=1}^4 \frac{1}{2}C_{ijk}[E_\alpha]_{jk} e_i \otimes E_\alpha\\
\label{eq:Cplus} & = \frac{1}{4}\sum_{\alpha=1}^3 \left(2W(E_\alpha)(\nabla f) + \frac{1}{3} E_\alpha \cdot (\Ric - R)(\nabla f) + \Ric \cdot E_\alpha (\nabla f)\right) \otimes E_\alpha.
\end{align}

\section{Smoothly diagonalizable Weyl curvature}
In view of the fact that $C^+ = -2\Div(W^+)$, we next derive an alternative expression of $C^+$ in terms of the eigenvalues $\lambda_\alpha$'s of $W^+$. Since differentiation will take place in the computation, we need to restrict the local frame $\{E_\alpha\}$ on the points in $\mathcal{O}$ at which $|\spec(W^+)|$ is locally constant. Some sketches of the following computations can be found in \cite{Der}. For reader's convenience, here we give more detail of the computations and revamp the index notations in a somewhat more systematic way.

Recall that
$$W^+ = \frac{1}{2}\sum_{\alpha=1}^3 \lambda_\alpha E_\alpha \otimes E_\alpha.$$

Using the notations defined in \eqref{eq:connection_E}, the covariant derivatives of $E_\alpha$'s are given by
\begin{equation*}
\nabla E_\alpha = \sum_{\beta=1}^3 \Omega_\alpha^\beta \otimes E_\beta.
\end{equation*}
Taking covariant derivative of $W^+$, we have
\begin{align*}
\nabla W^+ & = \frac{1}{2} \sum_{\alpha=1}^3 \left(d\lambda_\alpha \otimes E_\alpha \otimes E_\alpha + \sum_{\beta=1}^3\lambda_\alpha \Omega_\alpha^\beta \otimes (E_\alpha \otimes E_\beta + E_\beta \otimes E_\alpha)\right)\\
& = \frac{1}{2} \sum_{\alpha=1}^3 \left(d\lambda_\alpha \otimes E_\alpha \otimes E_\alpha + \sum_{\beta=1}^3\lambda_\alpha \Omega_\alpha^{\alpha+\beta} \otimes (E_\alpha \otimes E_{\alpha+\beta} + E_{\alpha+\beta} \otimes E_\alpha)\right).
\end{align*}
The last step follows from the cyclic permutation of the $\beta$ indices.

By \eqref{eq:switch_E} and applying cyclic permutation of indices (so as to produce a common factor $\otimes E_\alpha$ on each term), we have
\begin{align*}
& \sum_{\alpha=1}^3\sum_{\beta=1}^3\lambda_\alpha \Omega_\alpha^{\alpha+\beta} \otimes (E_\alpha \otimes E_{\alpha+\beta} + E_{\alpha+\beta} \otimes E_\alpha)\\
& = \sum_{\alpha=1}^3 \lambda_\alpha \Omega_{\alpha}^{\alpha+1} \otimes (E_\alpha \otimes E_{\alpha+1} + E_{\alpha+1} \otimes E_\alpha) \qquad\qquad (\beta = 1)\\
& + \sum_{\alpha=1}^3\lambda_\alpha \Omega_{\alpha}^{\alpha+2} \otimes (E_{\alpha} \otimes E_{\alpha+2} + E_{\alpha+2} \otimes E_{\alpha}) \qquad\qquad (\beta = 2)\\
& = \sum_{\alpha=1}^3 \lambda_{\alpha+2} \Omega_{\alpha+2}^{\alpha} \otimes E_{\alpha+2} \otimes E_\alpha + \sum_{\alpha=1}^3 \lambda_\alpha \Omega_{\alpha}^{\alpha+1} \otimes E_{\alpha+1} \otimes E_{\alpha} \qquad\qquad \text{(shifting indices)}\\
& + \sum_{\alpha=1}^3 \lambda_{\alpha+1}\Omega_{\alpha+1}^\alpha \otimes E_{\alpha+1} \otimes E_\alpha + \sum_{\alpha=1}^3 \lambda_{\alpha}\Omega_{\alpha}^{\alpha+2} \otimes E_{\alpha+2} \otimes E_{\alpha} \qquad\qquad \text{(shifting indices)}\\
& = \sum_{\alpha=1}^3 \{(\lambda_\alpha - \lambda_{\alpha+1}) \Omega_{\alpha}^{\alpha+1} \otimes E_{\alpha+1} + (\lambda_\alpha - \lambda_{\alpha+2}) \Omega_{\alpha}^{\alpha+2} \otimes E_{\alpha+2}\} \otimes E_{\alpha}.
\end{align*}
Therefore, we have
\begin{equation*}
\nabla W^+ = \frac{1}{2}\sum_{\alpha=1}^3 T_\alpha \otimes E_\alpha
\end{equation*}
where each $T_\alpha$ is a (0,3)-tensor defined by
$$T_\alpha = d\lambda_\alpha \otimes E_\alpha + (\lambda_\alpha - \lambda_{\alpha+1})\Omega_{\alpha}^{\alpha+1} \otimes E_{\alpha+1} + (\lambda_\alpha - \lambda_{\alpha+2})\Omega_{\alpha}^{\alpha+2} \otimes E_{\alpha+2}.$$

Let $[E_{\alpha}]_{ij}$'s, $[\Omega_{\alpha}^\beta]_p$'s and $[T_{\alpha}]_{ijk}$'s denote the local components of $E_\alpha$, $\Omega_\alpha^\beta$ and $T_\alpha$ respectively, i.e. 
$$E_\alpha = \sum_{i,j=1}^4 [E_{\alpha}]_{ij} e_i \otimes e_j, \quad \Omega_\alpha^\beta = \sum_{p=1}^4 [\Omega_{\alpha}^\beta]_p e_p, \quad T_\alpha = \sum_{i,j,k=1}^4 [T_{\alpha}]_{ijk} e_i \otimes e_j \otimes e_k.$$

In terms of these local components, $\nabla W^+$ can be expressed as
\begin{align*}
\nabla_p W^+_{ijkl} & = \left(\frac{1}{2}\sum_{\alpha=1}^3 T_\alpha \otimes E_\alpha\right)_{pijkl} = \frac{1}{2}\sum_{\alpha=1}^3 [T_{\alpha}]_{pij} [E_{\alpha}]_{kl}\\
& = \frac{1}{2}\sum_{\alpha=1}^3 e_p(\lambda_\alpha) [E_{\alpha}]_{ij} [E_{\alpha}]_{kl} + \frac{1}{2}\sum_{\alpha=1}^3 (\lambda_\alpha - \lambda_{\alpha+1})[\Omega_{\alpha}^{\alpha+1}]_p [E_{\alpha+1}]_{ij} [E_{\alpha}]_{kl}\\
& + \frac{1}{2}\sum_{\alpha=1}^3 (\lambda_\alpha - \lambda_{\alpha+2}) [\Omega_{\alpha}^{\alpha+2}]_p [E_{\alpha+2}]_{ij} [E_{\alpha}]_{kl}.
\end{align*}
Tracing $p$ and $i$, we obtain
\begin{align*}
\Div(W^+)_{jkl} & = \sum_{i=1}^4 \nabla_i W_{ijkl}^+\\
& = \frac{1}{2}\sum_{\alpha=1}^3 e_i(\lambda_\alpha) [E_{\alpha}]_{ij} [E_{\alpha}]_{kl} + \frac{1}{2}\sum_{\alpha=1}^3 (\lambda_\alpha - \lambda_{\alpha+1})[\Omega_{\alpha}^{\alpha+1}]_i [E_{\alpha+1}]_{ij} [E_{\alpha}]_{kl}\\
& + \frac{1}{2}\sum_{\alpha=1}^3 (\lambda_\alpha - \lambda_{\alpha+2}) [\Omega_{\alpha}^{\alpha+2}]_i [E_{\alpha+2}]_{ij} [E_{\alpha}]_{kl}\\
& = -\frac{1}{2}\sum_{\alpha=1}^3 [E_\alpha (d\lambda_\alpha)]_j  [E_{\alpha}]_{kl} - \frac{1}{2}\sum_{\alpha=1}^3 (\lambda_\alpha - \lambda_{\alpha+1}) [E_{\alpha+1} \Omega_{\alpha}^{\alpha+1}]_j [E_{\alpha}]_{kl} \\
& - \frac{1}{2}\sum_{\alpha=1}^3 (\lambda_\alpha - \lambda_{\alpha+2}) [E_{\alpha+2} \Omega_\alpha^{\alpha+2}]_j  [E_{\alpha}]_{kl}.
\end{align*}
Here $E_{\alpha+1}\Omega_{\alpha}^{\alpha+1}$ denotes the matrix product of $E_{\alpha+1}$ and $\Omega_{\alpha}^{\alpha+1}$.

Denote
$S_\alpha = \sum_{k=1}^4 [S_\alpha]_j e_j$ as a 1-form defined by:
$$[S_\alpha]_j = [E_\alpha (d\lambda_\alpha)]_k + (\lambda_{\alpha} - \lambda_{\alpha+1}) [E_{\alpha+1}\Omega^{\alpha+1}_\alpha]_j + (\lambda_{\alpha} - \lambda_{\alpha+2}) [E_{\alpha+2}\Omega^{\alpha+2}_\alpha]_j,$$
one get write $\Div(W^+)_{jkl} = -\frac{1}{2}\sum_{\alpha=1} [S_\alpha]_j [E_\alpha]_{kl}$.

By the quaternionic properties of $E_\alpha$'s, we have $E_{\alpha+1} = - E_\alpha E_{\alpha+2}$ and $E_{\alpha+2} = E_\alpha E_{\alpha+1}$. One can further rewrite $[S_\alpha]_k$ as:
$$[S_\alpha]_k = [E_\alpha U_\alpha]_k$$
where $U_\alpha = d\lambda_\alpha - (\lambda_{\alpha} - \lambda_{\alpha+1}) E_{\alpha+2}\Omega^{\alpha+1}_\alpha + (\lambda_{\alpha}-\lambda_{\alpha+2})E_{\alpha+1}\Omega^{\alpha+2}_\alpha$.

Therefore, one can express $\Div(W^+)$ locally as
\begin{equation}\label{eq:div_W+}
\Div(W^+) = -\frac{1}{2} \sum_{\alpha=1}^3 S_\alpha \otimes E_\alpha = - \frac{1}{2} \sum_{\alpha=1}^3 (E_\alpha U_\alpha) \otimes E_\alpha.
\end{equation}
Note that \eqref{eq:div_W+} is valid locally on an open set where $|\spec(W^+)|$ is constant and as such $\lambda_\alpha$'s are smooth.

\newpage
\section{Proof of Theorem \ref{thm:main}}
We call $W^+$ is \textit{degenerate} if $W^+$ has at most 2 distinct eigenvalues. One crucial consequence for both $g$ and $\Ric$ being $J$-invariant is that the self-dual Weyl tensor $W^+$ is everywhere degenerate, as established by Apostolov--Gauduchon in \cite{AG}.

For a Hermitian manifold $(M^4, g, J)$ in (real) dimension 4, the defect from being K\"ahler is measured by the \textit{Lee form} $\theta$ defined by
$$d\omega = \theta \wedge \omega.$$
Recall that $\omega$ is the Hermitian form defined by $\omega(X, Y) = g(X, JY)$.

When $W^+$ is degenerate, it was also established in \cite{AG} that $W^+$ has the following eigen-properties:
\begin{enumerate}
\item The eigenvalues of $W^+$ are given by: $\frac{\kappa}{6}$, $-\frac{\kappa}{12}$ and $-\frac{\kappa}{12}$, where $\kappa$ is the \textit{conformal scalar curvature} given by $\kappa = R - \frac{3}{2}(|\theta|^2 + 2\delta\theta)$. Here $\delta\theta$ is the co-differential of $\theta$.
\item The Hermitian form $\omega$ is an eigenform of $W^+$ with simple eigenvalue, i.e. $\frac{\kappa}{6}$
\end{enumerate}

Combining all these ingredients, we are now ready to give the proof of the main result.

\begin{proof}[Proof of Theorem \ref{thm:main}]
Let $\mathcal{O}$ be the open dense subset of $M$ consisting of points at which $|\spec(W^+)|$ is locally constant. By the $J$-invariant conditions on $g$ and $\Ric$, $|\spec(W^+)| \leq 2$ at every point on $M$.

Given any point $x \in \mathcal{O}$, we consider two separate cases: $|\spec(W_x^+)| = 1$ and $|\spec(W_x^+)| = 2$. We are going to show in both cases the desired identity $\Ric_0(\nabla f, \cdot) \wedge d\omega = 0$ holds. By continuity, it extends to the whole manifold $M$.

\textbf{Case 1:} $|\spec(W_x^+)| = 1$.

Let $\{e_i\}$ be a local orthonormal frame in a neighborhood of $x$ inside $\mathcal{O}$ and $E_\alpha$'s are local frame of $\wedge^+$ as defined in 
\eqref{eq:EF}. As $W^+$ is trace-free, the fact that $W_x^+$ only has one eigenvalue implies $W_x^+ = 0$ and hence $W^+ \equiv 0$ in a neighborhood of $x$.

When $W_x^+ = 0$, the relation between $\nabla f$ and $\Ric$ at $x$ was already discussed in works by Cao--Chen \cite{CaoChen} and Chen--Wang \cite{ChenWang} (and the references therein) for both shrinking and steady Ricci solitons. It was established in \cite{ChenWang} that if $\nabla f \not= 0$ at the point $x$, then $W^+_x = 0$ implies $\nabla f$ is an eigenvector of $\Ric$ (when considered as an endomorphism of $T_xM$). Furthermore, the Ricci tensor $\Ric$ at $x$ either equals $\frac{R}{4}I$ at the point $x$, or has two distinct eigenvalues with $1+3$ multiplities. However, the latter case is ruled out under our assumption that $\Ric$ is $J$-invariant, which implies the eigenvalues of $\Ric$ must come in pairs.

Therefore, whether or not $\nabla f = 0$ at $x \in \mathcal{O}$, we must have $\Ric_0(\nabla f, \cdot) = 0$ in this case.

We further remark that the above-mentioned eigen-properties of $\Ric$ can also be seen under the framework of this article yet in the same morale as in Chen-Wang's paper. First pick an orthonormal frame $\{e_i\}$ under which $\Ric$ at $x$ is diagonal. From \eqref{eq:Cplus}, we have
$$C^+ = \frac{1}{4} \sum_{\alpha=1}^3 \left(\frac{1}{3} E_\alpha \cdot (\Ric - R)(\nabla f) + \Ric \cdot E_\alpha (\nabla f)\right) \otimes E_\alpha.$$
Note that in this case $W^+ = 0$. On the other hand, since $C^+ = -2\Div(W^+) = 0$, it implies for any $\alpha$, we have
\begin{equation}\label{eq:W^+_zero}
\left(\frac{1}{3} E_\alpha \cdot (\Ric - R) + \Ric \cdot E_\alpha\right) (\nabla f) = 0.
\end{equation}
By expressing $\Ric_x$ in a diagonal form and rewriting \eqref{eq:W^+_zero} into a system of 12 equations with eigenvalues of $\Ric_x$ and components of $\nabla f$ as unknowns, one can recover $\Ric_0(\nabla f, \cdot) = 0$ at $x$ through elementary equation solvings.

\newpage
\textbf{Case 2:} $|\spec(W_x^+)| = 2$.

In this case, the eigenvalues of $W^+$ near $x$ are $\frac{\kappa}{6}$, $-\frac{\kappa}{12}$ and $-\frac{\kappa}{12}$, and the Hermitian form $\omega$ is an eigenform of $W^+$ with the simple eigenvalue $\frac{\kappa}{6}$. Without loss of generality, we let $e_2 = Je_1$ and $e_4 = Je_3$, then $\{e_1, Je_1, e_3, Je_3\}$ is a unitary frame near $x$ and  
$$E_1 = e_1 \wedge Je_1 + e_3 \wedge Je_3 =: \omega$$
is the Hermitian form. Since $E_2$ and $E_3$, defined according to \eqref{eq:EF}, are both orthogonal to $E_1$, they are eigenforms of $W^+$ with eigenvalues $-\frac{\kappa}{12}$.

Since $W^+$ is smoothly diagonalizable near $x$, the divergence $\Div(W^+)$ can be expressed as in \eqref{eq:div_W+}. By noting that $C^+ = -2\Div(W^+)$, we get for any $\alpha$:
\begin{align*}
& \frac{1}{4}\left\{2W(E_\alpha)(\nabla f) + \frac{1}{3} E_\alpha \cdot (\Ric - R)(\nabla f) + \Ric \cdot E_\alpha (\nabla f)\right\}\\
& = E_\alpha (d\lambda_\alpha - (\lambda_\alpha - \lambda_{\alpha+1}) E_{\alpha+2}\Omega_\alpha^{\alpha+1} + (\lambda_\alpha - \lambda_{\alpha+2})E_{\alpha+1}\Omega_{\alpha}^{\alpha+2})
\end{align*}
Multiplying $-4E_\alpha$ from the left and using the fact that $\lambda_1 = \frac{\kappa}{6}$, $\lambda_2 = \lambda_3 = -\frac{\kappa}{12}$, we get:
\begin{align}\label{eq:gradf_Omega}
\nonumber \left(\frac{\kappa}{3} I + \frac{1}{3} (\Ric - R) - \omega \cdot \Ric \cdot \omega \right) (\nabla f) & = \frac{2}{3}d\kappa - \kappa E_3 \Omega_1^2 + \kappa E_2 \Omega_1^3\\
\left(-\frac{\kappa}{6} I + \frac{1}{3} (\Ric - R) - E_2 \cdot \Ric \cdot E_2 \right) (\nabla f) & = -\frac{1}{3}d\kappa + \kappa E_3 \Omega^2_1\\
\nonumber \left(-\frac{\kappa}{6} I + \frac{1}{3}  (\Ric - R) - E_3 \cdot \Ric \cdot E_3 \right) (\nabla f) & = -\frac{1}{3}d\kappa - \kappa E_2 \Omega_1^3.
\end{align}
The complex structure $J$ acts on 2-forms $\beta$ by $J(\beta)(X,Y) = \beta(JX, Y)$. Using this definition, one can easily check $JE_2 = -E_3$.
Since $\Ric$ is $J$-invariant, we have $J^T \cdot \Ric \cdot J = \Ric$, or equivalently, $J \cdot \Ric \cdot J = - \Ric$. It is then straight-forward to check that
$$E_3 \cdot \Ric \cdot E_3 = JE_2 \cdot \Ric \cdot JE_2 = - JE_2 \cdot (J \cdot \Ric \cdot J) \cdot JE_2 = E_2 \cdot \Ric \cdot E_2.$$
Here we have used the fact that $JE_2J = E_2$ and $J^2 = -I$.

Thus from the second and third equations of \eqref{eq:gradf_Omega}, we get
$$E_3 \Omega_2^1 = - E_2\Omega_3^1.$$
Note also that $JE_1 = E_1J = I$, we also have 
$$\omega \cdot \Ric \cdot \omega = -\omega \cdot (J \cdot \Ric \cdot J) \cdot \omega = -\Ric,$$
the first equation of \eqref{eq:gradf_Omega} becomes
\begin{equation}\label{eq:R}
\left(\left(\frac{\kappa}{3} - \frac{1}{3}R\right) I + \frac{4}{3} \Ric \right) (\nabla f, \cdot) = \frac{2}{3}d\kappa - 2\kappa E_3 \Omega_1^2.
\end{equation}

Next we claim that $E_3 \Omega_1^2 = -\theta$, where $\theta$ is the Lee form. The proof goes by rewriting $d\omega = \theta \wedge \omega$ in terms of the local connection 1-forms $\Omega_\alpha^\beta$'s defined by \eqref{eq:connection_E}:

\begin{align*}
(d\omega)_{ijk} & = \nabla_i \omega_{jk} - \nabla_j \omega_{ik} + \nabla_k \omega_{ij}\\
& = [\Omega_1^2]_i[E_2]_{jk} + [\Omega_1^3]_i [E_3]_{jk}\\
& - [\Omega_1^2]_j[E_2]_{ik} - [\Omega_1^3]_j [E_3]_{ik}\\
& + [\Omega_1^2]_k[E_2]_{ij} + [\Omega_1^3]_k [E_3]_{ij},\\
(d\omega)_{ijk} [E_1]_{jk} & = [\Omega_1^2]_j [E_1 E_2]_{ji} + [\Omega_1^3]_j [E_1 E_3]_{ji}\\
& + [\Omega_1^2]_k [E_2 E_1]_{ik} + [\Omega_1^3]_k [E_3 E_1]_{ik}\\
& = -2[E_3 \Omega_1^2]_i + 2[E_2 \Omega_1^3]_i\\
& = -4[E_3 \Omega_1^2]_i.
\end{align*}
Note that we have used the fact that $E_1$ is orthogonal to $E_2$ and $E_3$ and hence $[E_1]_{jk} [E_2]_{jk} = -\tr(E_1 E_2)$ = 0. The same goes between $E_1$ and $E_3$.
On the other hand,
\begin{align*}
(\theta \wedge \omega)_{ijk} & = 2(\theta_i \omega_{jk} + \theta_j \omega_{ki} + \theta_k \omega_{ij}),\\
(\theta \wedge \omega)_{ijk} [E_1]_{jk} & = 2(\theta_i \omega_{jk} \omega_{jk} + \theta_j [\omega^2]_{ji} + \theta_k [\omega^2]_{ik})\\
& = 2(4\theta_i - \theta_j \delta_{ji} - \theta_k \delta_{ik}) = 4\theta_i.
\end{align*}
Equating $d\omega = \theta \wedge \omega$, we get $E_3 \Omega_1^2 = -\theta$.

Recall that the conformal scalar curvature $\kappa$ is given by
$$\kappa = R - \frac{3}{2}(|\theta|^2 + 2\delta\theta).$$
Combining this with \eqref{eq:R} and the fact that $2\Ric(\nabla f, \cdot) = dR$, we have
\begin{equation}\label{eq:kappa}
-\frac{1}{2} (|\theta|^2 + 2\delta\theta) df + d(|\theta|^2 + 2\delta\theta) = 2\kappa \theta.
\end{equation}
By differentiating both sides of \eqref{eq:kappa}, we have
$$-\frac{1}{2} d(|\theta|^2 + 2\delta\theta) \wedge df = 2d\kappa \wedge \theta + 2\kappa d\theta.$$
Combining this with \eqref{eq:kappa} and noting that $df \wedge df = 0$, we have $-\kappa \theta \wedge df = 2d\kappa \wedge \theta + 2\kappa d\theta$.
\begin{align}\label{eq:kappa_LHS}
-\kappa \theta \wedge df & = \left(R - \frac{3}{2}(|\theta|^2 + 2\delta\theta)\right) df \wedge \theta\\
\nonumber & = R df \wedge \theta - \frac{3}{2}(|\theta|^2 + 2\delta\theta) df \wedge \theta.
\end{align}
\begin{align}\label{eq:kappa_RHS}
2d\kappa \wedge \theta + 2\kappa d\theta & = (2dR - 3d(|\theta|^2 + 2\delta\theta)) \wedge \theta + 2\kappa d\theta\\
\nonumber & = 2 dR \wedge \theta - \frac{3}{2} (|\theta|^2 + 2\delta\theta) df \wedge \theta + 2\kappa d\theta.
\end{align}
The last equality of \eqref{eq:kappa_RHS} follows from taking wedge product on both sides of \eqref{eq:kappa} with $\theta$.

Equating \eqref{eq:kappa_LHS} and \eqref{eq:kappa_RHS}, we get $(Rdf - 2dR) \wedge \theta = 2\kappa d\theta$. Since $dR = 2\Ric(\nabla f, \cdot)$ and $df = g(\nabla f, \cdot)$, we get:
\begin{equation*}
\mathcal{\Ric}_0(\nabla f, \cdot) \wedge \theta = -\frac{1}{2}\kappa d\theta,
\end{equation*}
where $\mathcal{\Ric}_0$ denotes the trace-free Ricci tensor $\Ric - \frac{R}{4}g$.

Taking wedge product with $\omega$ from the right yields:
$$\Ric_0(\nabla f, \cdot) \wedge d\omega = -\frac{1}{2}\kappa d\theta \wedge \omega.$$

Finally, one can show $d\theta \wedge \omega = 0$ by the following argument:
\begin{align*}
d\omega & = \theta \wedge \omega\\
d^2\omega & = d\theta \wedge \omega - \theta \wedge \theta \wedge \omega\\
0 & = d\theta \wedge \omega.
\end{align*}
It proves the desired identity \eqref{eq:main} in this case.

The identity \eqref{eq:main} holds on the dense set $\mathcal{O}$. The proof of the theorem is completed by continuity.
\end{proof}
\begin{remark}
When $M$ is compact, the Lee form is closed, i.e. $d\theta = 0$. Therefore, one can obtain a stronger result than stated in the theorem:
$$\Ric_0(\nabla f, \cdot) \wedge \theta = 0.$$
\end{remark}

\bibliographystyle{amsplain}
\bibliography{citations}

\providecommand{\bysame}{\leavevmode\hbox to3em{\hrulefill}\thinspace}
\providecommand{\MR}{\relax\ifhmode\unskip\space\fi MR }
\providecommand{\MRhref}[2]{%
  \href{http://www.ams.org/mathscinet-getitem?mr=#1}{#2}
}
\providecommand{\href}[2]{#2}
\begin{thebibliography}{10}

\bibitem{AG}
V.~Apostolov and P.~Gauduchon, \emph{The {R}iemannian {G}oldberg-{S}achs
  theorem}, Internat. J. Math. \textbf{8} (1997), no.~4, 421--439. \MR{1460894
  (98g:53080)}

\bibitem{AHS}
M.~F. Atiyah, N.~J. Hitchin, and I.~M. Singer, \emph{Self-duality in
  four-dimensional {R}iemannian geometry}, Proc. Roy. Soc. London Ser. A
  \textbf{362} (1978), no.~1711, 425--461. \MR{506229 (80d:53023)}

\bibitem{CaoSoliton}
Huai-Dong Cao, \emph{Existence of gradient {K}{\"a}hler-{R}icci solitons},
  Elliptic and parabolic methods in geometry ({M}inneapolis, {MN}, 1994), A K
  Peters, Wellesley, MA, 1996, pp.~1--16. \MR{1417944 (98a:53058)}

\bibitem{CaoSurv}
\bysame, \emph{Recent progress on {R}icci solitons}, Recent advances in
  geometric analysis, Adv. Lect. Math. (ALM), vol.~11, Int. Press, Somerville,
  MA, 2010, pp.~1--38. \MR{2648937 (2011d:53061)}

\bibitem{CaoChen}
Huai-Dong Cao and Qiang Chen, \emph{On {B}ach-flat gradient shrinking {R}icci
  solitons}, Duke Math. J. \textbf{162} (2013), no.~6, 1149--1169. \MR{3053567}

\bibitem{CLW}
Xiuxiong Chen, Claude Lebrun, and Brian Weber, \emph{On conformally
  {K}{\"a}hler, {E}instein manifolds}, J. Amer. Math. Soc. \textbf{21} (2008),
  no.~4, 1137--1168. \MR{2425183 (2010h:53054)}

\bibitem{ChenWang}
Xiuxiong Chen and Yuanqi Wang, \emph{On four-dimensional anti-self-dual
  gradient ricci solitons},  (2011).

\bibitem{ChowBook}
Bennett Chow, Sun-Chin Chu, David Glickenstein, Christine Guenther, James
  Isenberg, Tom Ivey, Dan Knopf, Peng Lu, Feng Luo, and Lei Ni, \emph{The
  {R}icci flow: techniques and applications. {P}art {I}}, Mathematical Surveys
  and Monographs, vol. 135, American Mathematical Society, Providence, RI,
  2007, Geometric aspects. \MR{2302600 (2008f:53088)}

\bibitem{Der}
Andrzej Derdzi{{\'n}}ski, \emph{Self-dual {K}{\"a}hler manifolds and {E}instein
  manifolds of dimension four}, Compositio Math. \textbf{49} (1983), no.~3,
  405--433. \MR{707181 (84h:53060)}

\bibitem{FIK}
Mikhail Feldman, Tom Ilmanen, and Dan Knopf, \emph{Rotationally symmetric
  shrinking and expanding gradient {K}{\"a}hler-{R}icci solitons}, J.
  Differential Geom. \textbf{65} (2003), no.~2, 169--209. \MR{2058261
  (2005e:53102)}

\bibitem{KoisoSoliton}
Norihito Koiso, \emph{On rotationally symmetric {H}amilton's equation for
  {K}{\"a}hler-{E}instein metrics}, Recent topics in differential and analytic
  geometry, Adv. Stud. Pure Math., vol.~18, Academic Press, Boston, MA, 1990,
  pp.~327--337. \MR{1145263 (93d:53057)}

\bibitem{Page}
Don Page, \emph{A compact rotating gravitational instanton}, Physics Letters B
  \textbf{79} (1978), no.~3, 235--238.

\bibitem{WangZhu}
Xu-Jia Wang and Xiaohua Zhu, \emph{K{\"a}hler-{R}icci solitons on toric
  manifolds with positive first {C}hern class}, Adv. Math. \textbf{188} (2004),
  no.~1, 87--103. \MR{2084775 (2005d:53074)}

\end{thebibliography}
\end{document}